\documentclass[12pt,twoside]{article}

\usepackage[english]{babel}
\usepackage{amsmath}
\usepackage{amsfonts}
\usepackage{amssymb}
\usepackage{enumerate}
\usepackage{mathrsfs}
\usepackage{amssymb}
\usepackage[all]{xy}
\usepackage{graphics}

\usepackage{dsfont}
\usepackage{amsthm}

\date{}

\begin{document}

\newtheorem{Theorem}{\quad Theorem}[section]

\newtheorem{Definition}[Theorem]{\quad Definition}

\newtheorem{Corollary}[Theorem]{\quad Corollary}

\newtheorem{Lemma}[Theorem]{\quad Lemma}

\centerline{}

\centerline{}

\centerline {\Large{\bf A Note on a Conjecture about Commuting Graphs.}}

\centerline{}

\centerline{\bf {C. Miguel}}

\centerline{}

\centerline{ Instituto de Telecomunuca\c
c\~oes,}

\centerline{P\'olo de Covilh\~a,}

\centerline{celino@ubi.pt}

\begin{abstract} We prove that the diameter of the commuting graph of the full matrix ring over the real numbers is at most five. This answers, in the affirmative, a conjecture proposed by Akbari-Mohammadian-radjavi-Raja, for the special case of the field of real numbers.
\end{abstract}



\section{Introduction }For a ring $R$, we denote the center of $R$ by $Z(R)$, that is, $Z(R)=\{x\in R\;:\;xr=rx\;\forall r\in R\}$. The {\it commuting graph} of  $R$, denoted by $\Gamma(R)$, is a simple graph whose vertices are all non-central elements of $R$, and two distinct vertices $a$ and $b$ are adjacent if and only if $ab=ba$. In particular, the set of neighbors of a vertex $a$ is the set of all non-central elements of the centralizer of $a$ in $R$, that is, of $C_R(a)=\{x\in R : ax=xa\}$.  The commuting graph has been studied extensively in recent years by several authors, e.g.\cite{ak, ak1, ar,dol, do, do1, gi,ma, om}. Additional information about algebraic properties of the elements can be obtained by studying the properties of a commuting graph. For example, if   $R$ is a ring with identity such that $\Gamma(R)\cong\Gamma(M_2(\mathds F))$, for a finite field $\mathds F$, then $R\cong M_2(\mathds F)$, see \cite{ma}. It is conjectured that this is also true for the full matrix ring  $M_n(\mathds F)$, where $\mathds F$ is a finite field and $n>2$.

We next recall some concepts from graph theory. In a graph $G$, a {\it path}  is an ordered   sequence  $v_1-v_2-\ldots -v_l$ of distinct vertices of $G$  in which every two consecutive vertices are adjacent. The graph $G$ is called {\it connected} if there is at least one path between every pair of its vertices.
 The distance between two distinct vertices $u$ and $v$, denoted by $d(u , v)$, is the length of the shortest  path connecting them (if
such a path does not exists, then $d(u , v)=\infty$). The {\it diameter} of the graph is the longest distance between any two vertices of the graph $G$ and will be denoted by $diam(G)$.

 Much research has concerned the diameter of commuting graphs of certain classes of rings, e.g.\cite{ak2,dol, do,gi}. For the matrix ring $M_n(\mathds F)$, over an algebraically closed field $\mathds F$, the commuting graph is connected and the diameter is always equal to four, provided $n\geq 3$ \cite{ak2}. Note that for $n=2$ the commuting graph is always disconnected \cite[remark 8]{ak4}. If the field $\mathds F$ is not algebraically closed, the commuting graph $\Gamma(M_n(\mathds F))$ may be disconnected for an arbitrarily large integer $n$ \cite{ak3}. However, in $\cite{ak2}$ it is proved that for any field $\mathds F$ and $n\geq 3$, if $\Gamma(M_n(\mathds F))$ is connected, then the diameter is between four and six \cite{ak2}. Also in \cite{ak2} the authors advanced the following conjecture.

 \

{\bf Conjecture.}[Akbari, Mohammadian, Radjavi, Raja] Let $\mathds{F}$ be a field. If $\Gamma(M_n(\mathds F))$ is a connected graph, then its diameter is at most $5$.

\

The aim  of the present paper is to prove this conjecture for the special case of the field $\mathds R$ of real numbers. It was shown in \cite[remark 8]{ak3} that $\Gamma(M_n(\mathds R))$ is a connected graph, for any $n\geq 3$. It is well-know that  the central elements in the matrix ring $M_n(\mathds R)$ are the scalar matrices. Hence, the vertices in the graph $\Gamma(M_n(\mathds R))$ are the non-scalar matrices. Our main result is the following theorem.

\begin{Theorem}\label{t1}Let $\mathds R$ be the field of real numbers and $n\geq 3$. Then, $diam(\Gamma(M_n(\mathds R)))=4$ for $n\neq 4$ and $diam(\Gamma(M_4(\mathds R)))\leq 5$.
\end{Theorem}

We will show that Theorem \ref{t1} follows quickly from various results in the literature and from the real Jordan canonical form.

 We note that the Theorem \ref{t1} above generalizes easily to fields which have an algebraic closure which is a finite extension. This follows from the Artin-Schreier Theorem \cite[p.316]{ja} which asserts that such fields are precisely the real closed fields, which roughly speaking are the fields behaving like $\mathds R$ , and whose algebraic closures have degree $2$  and are given by adjoining a square root of $-1$. An example is the field of real algebraic numbers.

\section{Some Lemmas } In this section we assemble the tools that we require to prove Theorem\ref{t1}.
One of the key tools  is the real Jordan canonical form for a matrix over the field $\mathds R$  of real numbers. For the sake of completeness we recall it very briefly.
Let $\mathds F$ be an algebraically closed field and $A\in M_n(\mathds F)$. The well-known  Jordan canonical form theorem states that there is a regular matrix   $S\in M_n(\mathds F)$ such that
\begin{equation*}S^{-1}AS=diag(J_{n_1}(\lambda_1), \ldots , J_{n_k}(\lambda_k)).
\end{equation*}
Each $J_{n_i}(\lambda_i)$, $i=1, \ldots ,k$, is called a {\it Jordan block} of order $n_i$ and is of the form
\begin{equation*}J_{n_i}(\lambda_i)=\lambda_iI+N_i,
\end{equation*}
where $I$ is the identity matrix and each $N_i$, $i=1, \ldots , k$, is a square matrix whose only non-zero entries are $1's$ on the super-diagonal (i.e. just above the diagonal).
The matrix $S^{-1}AS=diag(J_{n_1}(\lambda_1), \ldots , J_{n_k}(\lambda_k))$ is called the {\it Jordan canonical form} of the matrix $A$.
 This canonical form was described by C. Jordan in 1870. If the field $\mathds F$ is not algebraically closed, the Jordan canonical form is no longer available for all matrices in $M_n(\mathds F)$.
 However, for matrices over the field of real numbers  the complex eigenvalues come in complex conjugate pairs, and this can be used to give a real Jordan canonical form for real matrices. Let $\lambda=a+ib$, where $b\neq 0$,  be a complex eigenvalue of a  matrix $A\in M_n(\mathds R)$.  Denote by $C(a , b)$ the $2$ by $2$ matrix
 \begin{equation*}C(a, b)=\left[\begin{array}{cc}a&b\\-b&a \end{array}\right]
 \end{equation*}
 and  by $C_k(a, b)$ the matrix of order $2k$ of the form
 \begin{equation*}C_k(a, b)=\left[\begin{array}{ccccc}C(a, b) &I  & 0 &
                           \ldots &0
                           \\ 0 & C(a, b)& I & \ddots & \vdots \\
                           0 & 0 & \ddots   & \ddots & 0 \\
                           \vdots & \ddots & \ddots & C(a, b) & I \\
                           0 & \ldots & 0 & 0& C(a, b)
                           \end{array}\right]
 \end{equation*}
 where $I$ is the identity matrix of order $2$. Then, there is a regular matrix  $S\in M_n(\mathds R)$ such that
 \begin{equation*}S^{-1}AS=diag(C_{n_1}(a_1, b_1), \ldots , C_{n_t}(a_t, b_t), J_{m_1}(\lambda_1), \ldots , J_{m_s}(\lambda_s))
 \end{equation*}
 where $a_j+ib_j$ are the nonreal eigenvalues of $A$, for $j=1, \ldots ,t$, and $\lambda_q$ are the real eigenvalues of $A$, for $q=1, \ldots ,s$. This is called the {\it real Jordan canonical form} of $A$. For more details about the real Jordan canonical form, including proofs, the reader is referred to section 6.7 of \cite{la}.

\

It was show in  \cite{dol}  that for  an algebraically closed field $\mathds F$, every matrix $A\in M_n(\mathds F)$ commutes with a rank one matrix. In the following result we show that for an  arbitrary field $\mathds F$, a matrix $A\in M_n(\mathds F)$ commutes with a rank one matrix if it has an eigenvalue in $\mathds F$. The proof is almost identical to the proof presented in \cite{dol} for algebraically closed fields.

\begin{Lemma}\label{l1} Let $\mathds F$ be a field and  $A\in M_n(\mathds F)$. If $A$ has an eigenvalue in $\mathds F$, then there exists a rank one matrix $X\in M_n(\mathds F)$ such that $d(A, X)\leq 1$.
\end{Lemma}
{\bf Proof.} Let $x$ and $y$ be eigenvectors of $A$ and $A^T$, respectively, corresponding to the same real eigenvalue $\lambda$. Then $X=xy^T$ is a rank one matrix with $AX=(Ax)y^T=\lambda xy^T=x(A^Ty)^T=XA$. $\Box$

\

Many of the results about the commuting graph of matrix rings over an algebraically closed field are obtained with the help of rank one matrices, see e.g. \cite{dol}. For a matrix $A\in M_n(\mathds R)$, if $A$ has no real eigenvalues, it is no longer true that $A$ commutes with a rank one matrix. Therefore, we can not apply to real matrices the same techniques that we apply in the case of matrices  over an  algebraically closed field. However, in the next lemma we show that every matrix $A\in M_n(\mathds R)$ commutes with a rank two matrix. This allows us  to adapt some techniques used with matrices over an algebraic closed field to matrices over the field of real numbers.

\begin{Lemma}\label{l2} Let $A\in M_n(\mathds R)$ be a matrix without real eigenvalues.  Then, there exists a rank two matrix $X\in M_n(\mathds R)$  such that  $d(A, X)\leq 1$.
\end{Lemma}
{\bf Proof.} Since the matrix $A$ has no real eigenvalues, it follows that $n$ is even and the real Jordan canonical form of $A$ is of the form
\begin{equation*}S^{-1}AS=diag(C_{n_1}(a_1, b_1), \ldots , C_{n_t}(a_t, b_t)).
\end{equation*}
Note that
 the field $\mathds C$ of complex numbers can be embedded   in the matrix ring $M_2(\mathds R)$ by the matrix-valued function $M:\mathds C\rightarrow M_2(\mathds R)$ defined by
\begin{equation}\label{e2}M(a+ib)=\left[\begin{array}{cc}a&b\\-b&a\end{array}\right].
\end{equation}
This embedding can be extended in the obvious way to an embedding function $\varphi : M_{n}(\mathds C)\rightarrow M_{2n}(\mathds R)$.

Now, observe that the real Jordan canonical form of $A$, that is, the matrix $S^{-1}AS$, belongs to the image of the embedding $\phi : M_{n/2}(\mathds C)\rightarrow M_{n}(\mathds R)$.  Let $E\in M_{n/2}(\mathds C)$ be such that $\varphi(E)=S^{-1}AS$. By Lemma \ref{l1} there exists a rank one matrix $R\in M_{n/2}(\mathds C)$ such that $d(E , R)\leq 1$. It follows that
$$ A(S\varphi(R)S^{-1})=(S\varphi(E)S^{-1})(S\varphi(R)S^{-1})=  S\varphi(E R)S^{-1}=(S\varphi(R)S^{-1})(S\varphi(E)S^{-1})=$$
$$=(S\varphi(R)S^{-1})A.
$$
To complete the proof we take $X=S\varphi(R)S^{-1}$.$\Box$

\

We conclude this section with three lemmas that have been proved in \cite{ak2} for the more general case of a division ring.

\begin{Lemma}\label{l3} Let $\mathds F$ be a field and $n\geq 2$. Suppose $A, B\in M_n(\mathds F)$ are two matrices such that $ker(A)\cap ker(B)\neq \{0\}$. Then, $C_{M_n(\mathds F)}(\{A, B\})$ contains at least one matrix with rank one.
\end{Lemma}

\begin{Lemma}\label{l4} Let $\mathds F$ be a field and $n\geq 3$. If $N, M\in M_n(\mathds F)$ are two non-zero matrices such that $M^2=N^2=0$, then $d(M , N)\leq 2$ in $\Gamma(M_n(\mathds F))$.\end{Lemma}
\begin{Lemma}\label{l5} Let $\mathds F$ be a field and $n\geq 3$. If $A, B\in M_n(\mathds F)$ are two non-scalar idempotent matrices, then
$d(A , B)\leq 2$ in $\Gamma(M_n(\mathds F))$. \end{Lemma}

\section{Proof of the Main Theorem}
In this section we prove the main theorem. Throughout this section $n\geq 3$ is a natural number.

\

{\bf Proof of Theorem \ref{t1}.} Suppose first $n\neq 4$. Let $A, B\in M_n(\mathds R)$ be matrices with no real eigenvalues. By Lemma \ref{l2} there exist matrices $X, Y\in M_n(\mathds R)$, with rank $2$ and such that $d(A, X)\leq 1$ and $d(B, Y)\leq 1$. Since $dim(Ker(A))=dim(Ker(B))=n-2$, it follows that $dim(Ker(X)\cap Ker(Y))\geq n-4$. Hence, $Ker(X)\cap Ker(Y)\neq \{0\}$. By Lemma \ref{l3} there is a matrix $Z\in C_{M_n(\mathds F)}(\{X, Y\})$, which implies that  $A-X-Z-Y-B$ is a path in $\Gamma(M_n(\mathds R))$.

Let now $A, B\in M_n(\mathds R)$ be matrices such that both have a real eigenvalue. By Lemma \ref{l1} there exist rank one matrices $X, Y\in M_n(\mathds R)$  such that $d(A , X)\leq 1$ and $d(B , Y)\leq 1$. Now we have $dim(Ker(X)\cap Ker(Y))\geq n-2$. Again by Lemma \ref{l3} there is a matrix $Z\in C_{M_n(\mathds R)}(\{X, Y\})$. Hence, $A-X-Z-Y-B$ is a path in $\Gamma(M_n(\mathds R))$.

Finally, let us assume that $A, B\in M_n(\mathds R)$ are such that $A$ has a real eigenvalue and $B$ has no real eigenvalues. By Lemma \ref{l1} there exists a  rank one matrix  $X\in M_n(\mathds R)$ such that $d(A , X)\leq 1$ and by Lemma \ref{l2} there exists a rank two matrix $Y\in M_n(\mathds R)$ such that $d(B , Y)\leq 1$. In this case we have $dim(Ker(X)\cap Ker(Y))\geq n-3$. Since $A$ has no real eigenvalues, it follows that $n\geq 4$.  Hence, $Ker(X)\cap Ker(Y)\neq\{ 0\}$. Again, there exists $Z\in M_n(\mathds R)$ such that $A-X-Z-Y-A$ is a path in $\Gamma(M_n(\mathds R))$. So, we have proved that for $n\neq 4$ $diam(\Gamma(M_n(\mathds R)))=4$.

 Suppose now that $n=4$ and let $A, B\in M_4(\mathds R)$. Assume first that both $A$ and $B$  have a real eigenvalue. As in the previous case, there exist rank one matrices $X, Y\in M_4(\mathds R)$  such that $d(A , X)\leq 1$ and $d(B , Y)\leq 1$. Since $dim(Ker(X)\cap Ker(Y))\geq 2$, it follows that there exists $Z\in M_4(\mathds R)$ such that $A-X-Z-Y-B$ is a path $\Gamma(M_4(\mathds R))$. If $A$ has a real eigenvalue and $B$ has no real eigenvalue, then there exist a rank one matrix $X\in M_4(\mathds R)$ and a rank two matrix $Y\in M_4(\mathds R)$ such that $d(A , X)\leq 1$ and $d(B , Y)\leq 1$. Since in this case we have $dim(Ker(A)\cap Ker(B))\geq 1$, it follows that there exists $Z\in M_4(\mathds R)$ such that $A-X-Z-Y-B$ is a path in $\Gamma(M_4(\mathds R))$.

 Finally, suppose that both $A$ and $B$ have no real eigenvalues.
 There are three possible cases  for the real Jordan form  of a matrix in $ M_4(\mathds R)$, namely

\begin{equation*}\left[\begin{array}{cc}R_1&0\\0&R_1\end{array}\right],\;\;\;\;\left[\begin{array}{cc}R_1&0\\0&R_2\end{array}\right]
\;\;or\;\;\left[\begin{array}{cc}R_1&I\\0&R_1\end{array}\right],
\end{equation*}
where $I$ is the identity matrix of order $2$ and $R_i$, for $i=1, 2$, is a square matrix of order 2 of the type $C(a, b)$. Therefore, for two matrices $A, B\in M_4(\mathds R)$ we have six possible cases. Let us study each case separately.

{\it Case $1$.}

\begin{equation*}A=S\left[\begin{array}{cc}A_1&0\\0&A_1\end{array}\right]S^{-1}\;\;\;and\;\;\;
B=L\left[\begin{array}{cc}B_1&0\\0&B_1\end{array}\right]L^{-1}.
    \end{equation*}
  In this case the matrices $A$ and $B$ commute with the following idempotent matrices
  \begin{equation*}I_1=S\left[\begin{array}{cc}I&0\\0&0\end{array}\right]S^{-1}\;\;\;and\;\;\;
I_2=L\left[\begin{array}{cc}I&0\\0&0\end{array}\right]L^{-1},
    \end{equation*}
    respectively. Since by Lemma \ref{l5} there is a non scalar matrix $X\in C_{M_4(\mathds R)}(\{I_1 , I_2\})$, it follows that $A-I_1-X-I_2-B$ is a path in $\Gamma(M_4(\mathds R))$.

    {\it Case $2$.}

    \begin{equation*}A=S\left[\begin{array}{cc}A_1&0\\0&A_1\end{array}\right]S^{-1}\;\;\;and\;\;\;
B=L\left[\begin{array}{cc}B_1&0\\0&B_2\end{array}\right]L^{-1}.
    \end{equation*}
    This case is similar to the case $1$.

     {\it Case $3$.}

     \begin{equation*}A=S\left[\begin{array}{cc}A_1&0\\0&A_1\end{array}\right]S^{-1}\;\;\;and\;\;\;
B=L\left[\begin{array}{cc}B_1&I\\0&B_1\end{array}\right]L^{-1}.
    \end{equation*}
    In this case the matrices $A$ and $B$ commute with the following nilpotent matrices
  \begin{equation*}N_1=S\left[\begin{array}{cc}0&I\\0&0\end{array}\right]S^{-1}\;\;\;and\;\;\;
N_2=L\left[\begin{array}{cc}0&I\\0&0\end{array}\right]L^{-1},
    \end{equation*}
    respectively. Since $N_1^2=N_2^2=0$, by Lemma \ref{l4}, there is a non scalar matrix $X\in C_{M_4(\mathds R)}(\{N_1 , N_2\})$. Hence, $A-N_1-X-N_2-B$ is a path in $\Gamma(M_4(\mathds R))$.

    {\it Case $4$.}

     \begin{equation*}A=S\left[\begin{array}{cc}A_1&0\\0&A_2\end{array}\right]S^{-1}\;\;\;and\;\;\;
B=L\left[\begin{array}{cc}B_1&0\\0&B_2\end{array}\right]L^{-1}.
    \end{equation*}
    This case is similar to the case $1$.

    {\it Case $5$.}

     \begin{equation*}A=S\left[\begin{array}{cc}A_1&0\\0&A_2\end{array}\right]S^{-1}\;\;\;and\;\;\;
B=L\left[\begin{array}{cc}B_1&I\\0&B_1\end{array}\right]L^{-1}.
    \end{equation*}
    In this case  $A$  commutes with the following idempotent matrix
  \begin{equation*}I_1=S\left[\begin{array}{cc}0&0\\0&I\end{array}\right]S^{-1}.
    \end{equation*}
    The matrix $B$ commutes with the following matrix

  \begin{equation*}B^\prime=L\left[\begin{array}{cc}B_1&0\\0&B_1\end{array}\right]L^{-1}.
    \end{equation*}
    Since $B^\prime$ commutes with the idempotent
  \begin{equation*}I_2=L\left[\begin{array}{cc}I&0\\0&0\end{array}\right]L^{-1},
    \end{equation*}
    it follows that there is a path $A-I_1-X-I_2-B^\prime-B$ in $\Gamma(M_4(\mathds R))$.

    {\it Case $6$.}

    \begin{equation*}A=S\left[\begin{array}{cc}A_1&I\\0&A_1\end{array}\right]S^{-1}\;\;\;and\;\;\;
B=L\left[\begin{array}{cc}B_1&I\\0&B_1\end{array}\right]L^{-1}.
    \end{equation*}
    This case is similar to the case $3$. The proof is completed.$\Box$

\

{\bf Remark:} We were unable to prove that the distance five, between two vertices in the graph $\Gamma(M_4(\mathds R))$,  is attained. Therefore, the question  $diam(\Gamma(M_4(\mathds R)))=5$ or  $diam(\Gamma(M_4(\mathds R)))=4$ remains  open.

 \bibliographystyle{plain}

\end{document}